\documentclass[11pt,reqno]{amsart}

\makeatletter
\@namedef{subjclassname@2020}{%
  \textup{2020} Mathematics Subject Classification}
\makeatother

 \usepackage{curves} 
\usepackage{graphicx} 
\usepackage{amsmath, amsfonts, amssymb, amscd,graphics, graphpap, tikz, color, xcolor, cancel,enumerate,enumitem, hyperref, mathabx, mathrsfs, mathtools}
\usepackage[scr=rsfs,cal=boondox]{mathalfa}

\usetikzlibrary{arrows.meta,arrows}
\usetikzlibrary{decorations.pathreplacing,decorations.markings}
\tikzset{arrow data/.style 2 args={%
      decoration={%
         markings,
         mark=at position #1 with \arrow{#2}},
         postaction=decorate}
      }%
      
\usepackage[mathscr]{euscript}
  \usepackage[normalem]{ulem}
  \usepackage{soul}

\usepackage{curves} 
\usepackage{graphicx}

\makeatletter
\newcommand{\doublewidetilde}[1]{{%
  \mathpalette\double@widetilde{#1}%
}}

\textwidth= 
148mm
\numberwithin{equation}{section}

\theoremstyle{plain}
\newtheorem{theo}{Theorem}[section]
\newtheorem{lem}[theo]{Lemma}
\newtheorem{prop}[theo]{Proposition}

\newtheorem{cor}[theo]{Corollary}

\theoremstyle{definition}
\newtheorem{rem}[theo]{Remark}

\newtheorem{example}[theo]{Example}
\newtheorem{definition}[theo]{Definition}

\newenvironment{pf}{\noindent{\it Proof.\,}}{\hfill $\square$\par \medskip}

\theoremstyle{plain}

\theoremstyle{definition}

\newcommand{\rank}{\operatorname{rank}}
\newcommand{\beq}{\begin{equation}}
\newcommand{\eeq}{\end{equation}}

\renewcommand{\d}{\delta}

\newcommand{\g}{\gamma}
\newcommand{\h}{\eta}

\renewcommand{\l}{\lambda}

\renewcommand{\r}{\rho}
\newcommand{\s}{\sigma}

\newcommand{\bR}{\mathbb{R}}

\newcommand{\bN}{\mathbb{N}}

\newcommand{\bV}{\mathbb{V}}

\newcommand{\gr}{\mathfrak{r}}

\newcommand{\cA}{\mathscr{A}}
\newcommand{\cB}{\mathcal B}
\newcommand{\cC}{\mathcal{C}}
\newcommand{\cD}{\mathscr{D}}

\newcommand{\cL}{\mathscr{L}}
\newcommand{\cM}{\mathscr{M}}

\newcommand{\cS}{\mathscr{S}}
\newcommand{\cT}{\mathscr{T}}
\newcommand{\cU}{\mathscr{U}}
\newcommand{\cV}{\mathscr{V}}
\newcommand{\cW}{\mathscr{W}}

\newcommand{\p}{\partial}

\renewcommand{\square}{\kern1pt\vbox
{\hrule height 0.6pt\hbox{\vrule width 0.6pt\hskip 3pt
\vbox{\vskip 6pt}\hskip 3pt\vrule width 0.6pt}\hrule height0.6pt}\kern1pt}

\DeclareMathOperator\Id{Id}

\renewcommand\Im{\operatorname{Im}}

\newcommand{\wt}{\widetilde}

\newcommand{\wc}{\widecheck}

\newcommand{\bt}{\begin{theo}\ \ }
\newcommand{\et}{\end{theo}}
\newcommand{\bp}{\begin{prop}\ \ }
\newcommand{\ep}{\end{prop}}
\newcommand{\bc}{\begin{cor}\ \ }
\newcommand{\ec}{\end{cor}}
\newcommand{\bl}{\begin{lem}\ \ }
\newcommand{\el}{\end{lem}}
\newcommand{\bd}{\begin{definition}}
\newcommand{\ed}{\end{definition}}
\newcommand{\be}{\begin{equation}}
\newcommand{\ee}{\end{equation}}

\def\<#1,#2>{\langle\,#1,\,#2\,\rangle}
\newcommand{\arr}{\begin{array}{rlll}}
\newcommand{\ea}{\end{array}}
\newcommand{\bea}{\begin{eqnarray}}
\newcommand{\eea}{\end{eqnarray}}
\newcommand{\bean}{\begin{eqnarray*}}
\newcommand{\eean}{\end{eqnarray*}}

\renewcommand{\=}{:=}

\newcommand{\ve}{\varepsilon}

  \newcommand{\vertiii}[1]{{\left\vert\kern-0.25ex\left\vert\kern-0.25ex\left\vert #1 
    \right\vert\kern-0.25ex\right\vert\kern-0.25ex\right\vert}}

\def\sideremark#1{\ifvmode\leavevmode\fi\vadjust{
\vbox to0pt{\hbox to 0pt{\hskip\hsize\hskip1em
\vbox{\hsize3cm\tiny\raggedright\pretolerance10000
\noindent #1\hfill}\hss}\vbox to8pt{\vfil}\vss}}}

\title[Proving the Chow-Rashevski\v\i Theorem \`a la Rashevski\u\i]{Proving the  Chow-Rashevski\u\i\  Theorem\\ \`a la Rashevski\u\i}
 \author[Cristina Giannotti,  Andrea Spiro and Marta Zoppello]{Cristina Giannotti \quad  Andrea Spiro \quad Marta Zoppello}
 
 
\begin{document}

\begin{abstract} 
We give a new independent proof of  a generalised version of the  theorem by Rashevski\u\i, which appeared in   [Uch. Zapiski Ped. Inst. K. {\bf 
2} 
(1938), 83 -- 94]  and from which the classical Chow-Rashevski\u\i\  Theorem follows as  a corollary. The proof is  structured  to allow generalisations to the case of  orbits of  compositions of flows in absence of  group structures, thus  appropriate   for   applications in Control Theory.   In fact,  the  same  structure of the proof  has been successfully exploited in [C. Giannotti, A. Spiro and M. Zoppello,   arXiv 2401.07555 \&  2401.07560 (2024)] to determine new controllability criteria for real analytic non-linear control systems. It also yields a corollary, which can be used to derive results under lower regularity assumptions, as it is illustrated by   a simple explicit example. 
\end{abstract}
\maketitle

\section{Introduction}
One of the most frequent  way to state  the Chow-Rashevski\u\i\  Theorem is the following: {\it If  $\cM$ is  a connected $n$-dimensional manifold    equipped with a  bracket generating distribution $\cD \subset T\cM$, 
 then any two points of $\cM$   can be  connected by a piecewise  smooth  $\cD$-path,  i.e. by a path whose tangent vectors are all   in  $\cD$.} However, in the original papers \cite{Ra, Ch},  Rashevski\u\i\  and Chow proved  two different theorems  on the  set  $\cS^{(x_o)}$ of the points of a  manifold with distribution $(\cM, \cD)$,  that    can be joined to a given $x_o \in \cM$ by means of $\cD$-paths and the above  statement of Chow-Rashevski\u\i\  Theorem  is  a corollary of such  two different results.   Curiously enough,  the two    theorems   have different hypotheses and are equivalent one to the other only if a certain  additional condition holds.  \par
 \smallskip
 On the one hand,  Rashevski\u\i's Theorem has the hypothesis  that  the  enlarged distribution $\cD^{(\text{Lie})} $,  generated by all vector fields  in $\cD$   and  all possible  iterated   Lie brackets between them,  has constant rank
 (\footnote{Strictly speaking,  in \cite{Ra} the author assumes  that  $\rank \cD^{(\text{Lie})}_x = \dim \cM$ at all points, but its proof works perfectly  well under the   weaker   assumption that  $\rank \cD^{(\text{Lie})}_x$  is constant.}).   In  this case,   he proves that, for any $x_o \in \cM$,   the corresponding   set $\cS^{(x_o)}$ is   the maximal integral submanifold of $\cD^{(\text{Lie})}$  (which  is easily seen to be involutive) that passes through $x_o$.   On the other hand,  the hypothesis in   Chow's Theorem is  that  a different enlarged  distribution has constant rank,  i.e. the distribution  $\cD^{(\text{Flows})}$ spanned by all vector fields  in $\cD$ together  with the vector fields of the form $\Phi^X_{t*}(Y)$, where $X, Y$ are  in $ \cD$ and $ \Phi^X_{t*}(Y)$ is the pushed forward of $Y$ by the flow of $X$ for all sufficiently small  $t$'s.  Under this assumption,  Chow proves that $\cD^{(\text{Flows})}$ is involutive and that,  if $\cD^{\text{Flows}}$ has constant rank, for any $x_o \in \cM$,  the corresponding set  $\cS^{(x_o)}$ coincides with the maximal   integral submanifold  of $\cD^{(\text{Flows})}$ passing through $x_o$. On this regard, it  should be observed that  the celebrated   Orbit Theorem by Sussmann \cite{Su} (proved several years later than Rashevski\u\i\  and Chow's papers)  shows that the constant rank condition on  $\cD^{(\text{Flows})}$ can be removed  and  that   $\cS^{(x_o)}$ is {\it always} a  maximal   integral submanifold  for the (possibly non-regular)  distribution $\cD^{(\text{Flows})}$.  \par
 Since the distribution  $\cD^{(\text{Lie})}$ is  always a  subdistribution of   $\cD^{(\text{Flows})}$,  by Rashevski\u\i's and Chow and Sussmann's results,  whenever  $\cD^{(\text{Lie})}$ has a constant rank, say $p$ --  so that  $\cS^{(x_o)}$ is   a $p$-dimensional maximal  integral submanifold for both  $\cD^{(\text{Lie})}$ and   $\cD^{(\text{Flows})}$  -- 
  then also $\cD^{(\text{Flows})}$  has constant rank  and its rank is  $p$ as well. Thus, in  these cases,   Rashevski\u\i's and Chow's theorems are  equivalent. Differences occur only  when  $\cD^{(\text{Lie})}$ has non-constant rank. In these situations  only Chow's (or Sussmann's) theorem applies, meaning that,  in a sense, the result by Chow  is  more powerful than the one by Rashevski\u\i. On the other hand,  the dimensions of   the spaces  $\cD^{(\text{Lie})}|_x$, $x \in \cM$,  (and hence the dimension of the integral submanifold $\cS^{(x_o)}$ when  $\cD^{(\text{Lie})}$  is regular)  is undoubtedly much easier to be determined than the dimension of $\cD^{(\text{Flow})}|_x$:  For the first one needs just  to compute Lie brackets,   while for the  second  it is necessary   to integrate the vector fields in $\cD$,  determine their flows 
 and, finally,  compute  the pushed-forward vector fields   $\Phi^X_{t*}(Y)$. Due to this, wherever  Rashevski\u\i's condition holds  (note that, by    semicontinuity of  the rank function,  for any point  of  an open and dense subset of  $(\cM, \cD)$, there is a neighbourhood $\cU$ of that point  where the restricted distribution   $\cD^{(\text{Lie})}|_\cU$  has constant rank), Rashevski\u\i's theorem is much  handier than Chow's. \par
 \smallskip
 In this paper, we present a proof of (a minor generalisation of)  Rashevski\u\i's theorem, which we constructed independently,  before we were  able to  read  \cite{Ra}
 (\footnote{It was quite hard   to get a  copy of \cite{Ra}. We succeeded in getting a copy   from the Russian State Library with the help of a  long chain of  friends.}). In fact,  our proof is  different from Rashevski\u\i's one and  is crucially rooted on a  lemma on 
 iterated Lie brackets, which was also  proved  by Feleqi and Rampazzo in \cite{FR} and  has an independent interest on its own (\footnote{Also Rashevski\u\i\  built his proof on  a   lemma on iterated Lie brackets. His lemma  is surely different from ours,   but  it is very close to it in the spirit.}).  
The circle of ideas we follow  is  designed to    allow several kinds of  modifications  and  generalisations, more appropriate  for studies of     orbits of points  under families of  compositions of flows which  do not have a natural  group structure, as it often occurs in Control Theory. Modifications and generalisation  of this kind have been  successfully exploited in  \cite{GSZ1, GSZ2} (see also \cite{BGSZ}).  As a corollary of our proof, we also get an estimate for the  radius of a ball, made of  points joinable  to its centre through $\cD$-paths. The estimate is designed to be applicable  in discussions of  deformations of distributions and  of Hausdorff limits of  orbits.  At the end of the last section,   we give a simple explicit example which   illustrates  one of the ways this   estimate might be exploited. 
\par
\smallskip
The structure of the paper is as follows. After a preliminary section \S \ref{sect2}, where we review the notion of piecewise smooth  curves tangent to a distribution and we fix   once and for all the  definition of smooth  non-regular distribution, in \S \ref{derived}  we introduce the notion of distribution of uniform $\mu$-type and   state Rashevski\u\i's Theorem.  
In \S \ref{appendix} and \S \ref{corchow}, the announced proof  and its corollary (Corollary \ref{cor43}) on the radii of the balls, made  of points that are joinable to the centre via $\cD$-paths, are given.\par

\medskip
\section{Curves tangent to distributions} \label{sect2}
\subsection{Piecewise regular  curves} Let  $\cM$ be a connected $N$-dimensional manifold of class $\cC^\infty$.  In this paper by  {\it (regular) curve}  we always mean a subset of $\cM$,  which is   homeomorphic to a closed bounded interval of $\bR$ and equal to the image $\g([a, b])$  of a   smooth and regular  map $\g: [a, b] \subset \bR \to \cM$   (i.e.  a $\cC^1$ map with $\dot \g(t) \neq 0$ at all points). In other words, with the term ``curve'' we mean   the trace of a parameterised curve $\g(t)$ of class $\cC^1$, which has  nowhere vanishing velocity $\dot \g(t)$ (this is what is called   ``regular parameterisation'' throughout the paper) and which is  a homeomorphism between   its domain $[a, b] \subset \bR$ and its image $\g([a,b])$. 
Let $\g_1(t)$, $\g_2(s)$, $t \in [a, b]$,  $s \in [c,d]$, be  two  regular parameterisation of two curves (which, for simplicity,  we just denote by $\g_1$ and $\g_2$), such that the  endpoint   $\g_1(b)$ of the first curve is equal to the 
endpoint  $\g_2(c)$ of the second curve.  The  {\it  composition} $\g_1\ast \g_2$ is the  union of these two  curves. The curves  $\g_1$, $\g_2$ are the {\it regular arcs} of $\g_1 \ast \g_2$.   In a similar way  we  define the  {\it  composition} of a finite number of  regular parametrisations  $\g_1$, $\g_2$, $\ldots$, $\g_r$,  each of them sharing  one of its  endpoints with an endpoint of the succeeding one. A connected subset of $\cM$, which is obtained as composition of a finite collection of curves,  is called {\it piecewise regular curve}.\par
\smallskip
\subsection{Regular and singular distributions and their tangent curves}  
 A {\it regular distribution of rank $p$}  on a manifold $\cM$  is a smooth family $\cD$  of subspaces $\cD_x \subset T_x \cM$ of the tangent spaces of $\cM$ of constant dimension $p$. Here, by ``smooth family of subspaces'' we   mean a set of  subspaces  with the following property:   for any  point $x_o \in \cM$ 
there exists a neighbourhood $\cU$ and a  set of $\cC^\infty$ pointwise linearly independent vector fields $X_1, \ldots, X_p$ on $\cU$,  such that 
$$\cD_x = \langle X_1|_x, \ldots, X_p|_x \rangle\qquad \text{for any}\ x \in \cU\ .$$
Given a regular  distribution $\cD \subset T\cM$,  we will  use the short notation ``$X \in \cD$'' to indicate that the (local) vector field $X$ takes values in the vector spaces  $\cD_x$  of the distribution. \par
\smallskip
Generalisations of the notion of ``regular distribution'',  in which the  condition $\dim \cD_x = \text{constant}$  is not assumed, are possible,  but  demand some  care.  In this paper we   adopt  the  following definition,  which is a  variant of   those  considered in   \cite{Na, SJ, Su}. 
\begin{definition}\label{quasi-regular-def}  A {\it  quasi-regular set of   $\cC^\infty$  vector fields of rank $p$\/} on $\cM$ is a set $V$ of local vector fields of the following kind. There exist  
  an open cover  $\{\cU_A\}_{A \in J}$  of $\cM$ and a family    $\{ (X^{(A)}_{1}, \ldots, X^{(A)}_{p_A})\}_{A \in J}$  of tuples of cardinalities $p_A \geq p$,   of $\cC^\infty$  vector fields  -- one tuple for each   open set $\cU_A$ --  each of them containing a $p$-tuple, made of vector fields that are 
 pointwise linearly independent on some open and dense subset
of $\cU_A$,  and  satisfying the following  conditions: 
\begin{itemize}[leftmargin = 10pt]
\item  if $\cU_A \cap \cU_B \neq \emptyset$,  then  for any $x \in \cU_A \cap \cU_B$   one has $ X^{(A)}_{i}|_x = \cA^{(AB)}{}_{i}^j\big|_x X^{(B)}_{j}\big|_x$ for a $\cC^\infty$  matrix valued map
$\cA^{(AB)} : \cU_{A} \cap \cU_B \to \bR_{p_A \times p_B}$
\item the vector fields $Y \in V$  are exactly the local vector fields,   for which  any  restriction  $Y|_{\cU \cap \cU_A}$ to the intersection  between the   domain $\cU$  and a set  in  $\{\cU_A\}_{A \in J}$, has   the form
\beq Y|_{\cU_A \cap \cU} = Y^{(A)i} X^{(A)}_{i}\eeq
for  some $\cC^\infty$ functions $Y^{(A)i}$.  
\end{itemize}
The tuples $(X^{(A)}_{1}, \ldots, X^{(A)}_{p_A})$  are   called {\it sets of local  generators for $V$}.\par 
 A 
  {\it  generalised distribution of rank $p$}   is a pair $(V, \cD^V)$ given by  (i) a quasi-regular set of vector fields    $V$ of rank $p$  and  
(ii)  the uniquely associated family $\cD^V$ of   subspaces of tangent spaces of $\cM$,  given  by 
  $$\cD^V_x  = \{X_x\ ,\ X \in V \}\ ,\qquad x \in \cM\ .$$
 If $\dim \cD^V_x $  is constant, the   generalised distribution $(V, \cD)$ is called {\it regular}, otherwise it is called {\it singular}.
\end{definition}
 Notice that if  $(V, \cD^V)$ is   regular, then $V$ coincides with the set of local vector fields with  values in  the regular distribution $\cD = \cD^V$ and $(V, \cD^V)$ is  fully determined by  $\cD = \cD^V$.  On the other hand,  if  $(V, \cD^V)$ is   singular, there might be several different   generalised distributions,  all  of them  determining the same  family of tangent subspaces $\cD^V$. For instance
 the quasi-regular sets of local vector fields on $\bR^2$ 
 \begin{align*}
 &V = \Big\{\ X = f^1(x) \frac{\p}{\p x^1}  + f^2(x) x^1 \frac{\p}{\p x^2}\ ,\ f^i \ \text{smooth}\ \Big\}\ ,\\
 & \wt V = \Big\{\ X = g^1(x) \frac{\p}{\p x^1}  + g^2(x) (x^1)^2 \frac{\p}{\p x^2}\ ,\ g^i \ \text{smooth}\ \Big\}\ 
 \end{align*}
 determine the same families of subspaces $\cD^V = \cD^{\wt V}$, but $(V, \cD^V) \neq (\wt V, \cD^{\wt V})$ because    $\wt V$ is a proper subset of   $V$. 
\par

A  generalised distribution $(V, \cD^V)$   is called {\it non-integrable} if there exists at least one  pair of vector fields $X$, $Y \in V$, whose  Lie bracket $[X, Y]$  is not  in $V$.  It is called {\it integrable} or {\it involutive}, otherwise.  In case $(V, \cD^V)$ is regular and thus identifiable with the regular distribution $\cD \= \cD^V$, the non-integrability is  equivalent to  the existence of $X, Y \in \cD$ such that $[X, Y]_x \notin \cD|_x$ at some $x \in \cM$. \par
\smallskip
\begin{rem} The notion of generalised distributions  can be equivalently  defined   using the language of sheaves  as follows. First of all,  note that  if $(V, \cD^V)$ is a generalised distribution on  a manifold $\cM$,  then  the  collection of the sets of vector fields  in $V$ that are defined on open sets and the maps between pairs of such  collections,  determined by the usual  restriction map, give the ordered pair 
$$\bigg(\big\{ V|_{\cU}\big\}_{\cU\ \text{open}}\ , \big\{\r_{\cU, \cV}: V|_{\cU} \to V|_{\cV}\big\}_{\cV \subset \cU\ \text{open}}\bigg)\ , $$
which is a pre-sheaf and  corresponds to  a  subsheaf  $\pi: \bV^{(\cV, \cD)} \to \cM$ of the sheaf  $\pi: \cT \cM \to \cM$ of the  germs of  the  local $\cC^\infty$-vector fields of $\cM$. This subsheaf  satisfies the following two conditions: 
(a) there exists  an open dense  subset $\wt \cM \subset \cM$,  on which  each fiber $\bV|_x$, $x \in \wt \cM$,  span a subspace $\cD_x \subset T_x M$ of dimension $p$;  
  (b)  for any  $x \in \cM$ there is  a neighbourhood $\cU$ and a tuple  $(X_1, \ldots, X_{\bar p})$, $\bar p \geq p$,  of 
sections  of $\bV|_{\cU}$   generating   all other sections of  $\bV|_{\cU}$   over the ring of $\cC^\infty$ functions, i.e. the sheaf is  {\it locally finite generated}. Conversely, one can directly check that for any  subsheaf $\bV \subset \cT \cM$ satisfying (a) and (b),  the  set $V$ of all local sections  of $\bV$ together with  the  family of  the  vector spaces $\cD^V_x  =\text{Span}( \bV_x) \subset T_x M$ gives  a generalised distribution $(V, \cD^V)$   in the above sense.
\end{rem}

A  curve $\g$  is said to be {\it tangent to the generalised distribution $(V,\cD =  \cD^V)$}  if for one (hence,  for all) regular parameterisation $\g(t)$ of  the curve, the velocities $\dot \g(t)$  are such that  $\dot \g(t) \in \cD_{\g(t)}$  for any $t$.  A piecewise regular curve $\g = \g_1 \ast \g_2 \ast \ldots \ast \g_r$  is    {\it tangent to $(V, \cD = \cD^V)$} if each of its regular arcs has this property. We  call  any such $\g$  a {\it  $\cD$-path}.\par
The following relation between points of $\cM$
\beq\label{relation}  x \sim x' \  \text{if and only if} \  
\text{ there is a $\cD$-path  with endpoints  $x$ and $x'$}\eeq 
is an equivalence relation. Indeed,  in order to check that  \eqref{relation}  is   reflexive, symmetric and transitive, it suffices to recall  that each  (regular) curve,  which is tangent to  $\cD$ at all points,  is  simultaneously the image of a regular parameterisation $\g: [a, b] \to \cM$   and of its   ``inverse'' regular parameterisation  $\overline \g: [-b,-a] \to \cM$,   $\overline\g(t) \= \g(-t)$. The equivalence classes of such relation 
are called {\it $\cD$-path components of $\cM$}. If $\cM$ consists of just one $\cD$-path component, we say that $\cM$ is {\it $\cD$-path connected}.
\par
\medskip
\section{Distributions of uniform  $\mu$-type and Rashevski\u\i's Theorem} \label{derived}
\subsection{Distributions of uniform  $\mu$-type} \label{sect3.1}
Let $(V, \cD = \cD^V)$ be a generalised  distribution on a manifold $\cM$. We denote by    $V^{\operatorname{Lie}}$  the smallest family of   smooth local vector fields on $\cM$, which contains  all   local vector fields $X \in  V$ and  is closed under Lie brackets and under the standard operations of sum and multiplication by smooth real functions.  An equivalent {\it constructive} (and more formal) definition   of    $V^{\operatorname{Lie}}$ is  given   in   the language of ``germs''  and ``sheaves'' as follows.  \par
\smallskip
First of all,  let us fix   some notational convention. We denote sheaves of germs of smooth vector fields on $\cM$ by
   underlined calligraphic symbols,  as e.g. $\underline \cA$, $\underline \cB$, etc.  Moreover,  given  a generalised distribution $(V, \cD= \cD^V)$,  we denote by    $\underline{\cD}$  the sheaf of  the germs of the (local) smooth vector fields   that are in  $V$.  
We recall that  the  standard operations  $+$ and  $[\cdot, \cdot]$ of sum and Lie brackets between vector fields induce natural corresponding operations between  germs of vector fields. Given two sheaves $\underline \cA$,  $\underline \cB$ of germs of vector fields,  we denote by 
$$\underline \cA + \underline \cB\qquad \text{and}\qquad [\underline \cA, \underline \cB]$$
the sheaves given by    sums and  Lie brackets, respectively,   of the germs   in  $\underline \cA$ and $\underline \cB$.  In this way, 
for any generalised distribution $(V, \cD = \cD^V)$ on $ \cM$ we may consider the  following   increasing sequence of sheaves:
 \beq \label{41}
 \begin{split}
 &  \underline \cD^{(-1)}\= \underline \cD\ ,\\
 &\underline \cD^{(-2)} \=  \underline \cD^{(-1)} + [\underline \cD, \underline \cD^{(-1)}] \ ,\\
 &\underline \cD^{(-3)} \=   \underline \cD^{(-2)} + [\underline \cD, \underline \cD^{(-2)}]\ ,\\
 & \text{etc.}
 \end{split}
 \eeq
 A simple inductive argument based on the Jacobi identity shows that if $X, Y$ are two local smooth vector fields defined on a common  open set $\cU \subset \cM$ and with  germs  in $\underline \cD^{(-\ell)}$ and $\underline \cD^{(-m)}$, respectively, then their  Lie bracket $[X, Y]$ is  a vector field with germs in $\underline \cD^{(- \ell - m)}$. This means that 
 {\it the above described  class  $V^{\operatorname{Lie}}$  of local vector fields on $\cM$ coincide with the  class of all possible   local vector fields, whose germs  belong to one of the sheaves  $\underline \cD^{(-\ell)}$, $\ell \geq 1$}. In particular, $V^{\operatorname{Lie}}$ is closed under sum and Lie brackets of pairs of vector fields defined on common open sets.\par
\medskip
 We now observe  that, for any fixed point $x \in \cM$, the values at $x$ of the germs of  $\underline \cD^{(-\ell)}$ give a linear subspace $\underline \cD^{(-\ell)}|_x$ of $T_x \cM$, and determine a nested  sequence of linear subspaces 
  \beq \label{42}\underline \cD_x = \underline \cD^{(-1)}|_x \subset  \underline \cD^{(-2)}|_x \subset \underline \cD^{(-3)}|_x \subset \underline \cD^{(-4)}|_x \subset \ldots \ .
 \eeq
 We remark that if $\ell$ is such that $\underline \cD^{(-\ell)}|_x = \underline \cD^{(-\ell-1)}|_x$  at all points $x \in \cM$, then $\underline \cD^{(-\ell)} =  \underline \cD^{(-\ell-1)}$ 
  and   
  $\underline \cD^{(-\ell)} = \underline \cD^{(-\ell-s)}$ for any $s \geq 1$.\par
If there is a finite number   $1 \leq \nu $ for which the dimensions of  the spaces $ \underline \cD^{(-\ell)}|_x$  stabilise for all $x \in  \cM$ and $\ell \geq \nu$, we call     {\it  minimal depth} of $\cD$  the smallest integer  $\mu$  with this property.   In the other cases, we say that the {\it minimal depth is infinite}.\par
 \smallskip
 \begin{definition} \label{defuniform} A generalised distribution $(V, \cD = \cD^V)$   is called {\it of uniform  $\mu$-type} if  
 \begin{itemize}[leftmargin = 15pt]
 \item it is of finite minimal depth $\mu$; 
 \item  the family of vector spaces $(\cD_{-\mu})_x \= \underline \cD^{(-\mu)}\big|_x $, $x \in \cM$, is  a  regular distribution.  
 \end{itemize} 
 If  $(V, \cD = \cD^V)$ is of uniform $\mu$-type and   $\cD_{-\mu}|_x = T_x\cM$ for  any $x$, the generalised distribution $(V, \cD = \cD^V)$  is   called  {\it bracket generating}. \par
\end{definition}
Notice that   a generalised  distribution $(V, \cD = \cD^V)$ is of uniform $\mu$-type  of  depth  $\mu = 1$ if and only if $\cD$ is regular and  involutive. In  the other cases (regardless of being of uniform $\mu$-type or not) it is  non-integrable. \par
\smallskip
The following simple lemma can be directly checked. 
\begin{lem}  \label{lemmino} If $(V, \cD = \cD^V)$ is a generalised distribution  of uniform $\mu$-type, for any $x \in \cM$, there exists a neighbourhood $\cU \subset \cM$ and a set  of local generators $(Y_1, \ldots, Y_p)$ for $V|_{\cU}$, such that 
the set of all vector fields of the form
$$Y_{i_1}\ ,\qquad [Y_{i_1}, Y_{i_2}]\ ,\qquad \ldots \qquad [Y_{i_1}, [Y_{i_2}, [\ldots [Y_{i_{\mu-1}}, Y_{i_\mu}]\ldots ]]]\ ,\qquad 1 \leq  i_\ell \leq p\ ,$$
constitute  a set of local generators (not necessarily pointwise linearly independent) for the regular distribution $\cD_{-\mu}$.
\end{lem}
\par
\smallskip
\subsection{Rashevski\u\i's Theorem  for generalised distributions of uniform  $\mu$-type} \label{statementChow}
Let $(V, \cD)$ be a generalised distribution on $\cM$ of uniform $\mu$-type. In this case, the distribution $\cD_{-\mu}$ coincides with the distribution $\cD^{(\operatorname{Lie})}$, described in the Introduction and  the class of  the local vector fields in $\cD_{-\mu}$  is  the above considered set $V^{\operatorname{Lie}}$. In particular, {\it   $\cD_{-\mu}$ is  involutive}.   Being also  regular, by   the  classical Frobenius Theorem,  for any  $y_o \in \cM$  there exists a unique  maximal  integral leaf  $ \cS^{(y_o)}$  of   $\cD_{-\mu}$ passing through $y_o$ and   such a leaf coincides  with the $\cD_{-\mu}$-path component of $y_o$ (see e.g. \cite{Wa}). \par
Moreover, if we denote by $M$  the rank of $\cD_{-\mu}$, $M \leq \dim \cM$, any such integral leaf $\cS^{(y_o)} $ has  the form 
$\cS^{(y_o)} = \imath(\wt \cS)$ for  an injective immersion $\imath: \wt \cS \to \cM$  of  an $M$-dimensional  manifold $\wt \cS$. In particular, for any $y \in \cS^{(y_o)}$,   there is  a neighbourhood $\wt \cV \subset \wt \cS$ of  the point $\wt y = \imath^{-1}(y)$ such that   $ \imath(\wt \cV) $  is an $M$-dimensional  embedded submanifold  passing through $y$.  In what follows, any set $\imath(\wt \cV)$ determined in  this way will be  called {\it  neighbourhood of $y$  of  the intrinsic topology of $\cS^{(y_o)}$}. 
\par
\medskip
We are now  ready to state our   generalisation  of Rashevski\u\i's  result, originally established just for regular distributions.  
\begin{theo}
\label{Chow} If  $(V, \cD = \cD^V)$ is  a generalised distribution   of  uniform  $\mu$-type, 
 the $\cD$-path  components of $\cM$ are equal to  the $\cD_{-\mu}$-path components, that is to the  maximal integral leaves of the involutive   distribution $\cD_{-\mu}$. 
In particular, $\cM$ is $\cD$-path connected if and only if $\cD$ is bracket generating. 
\end{theo}
\par
\medskip
\section{Proof of  Theorem \ref{Chow}}
\label{appendix}
In this section  we give a   proof of Theorem \ref{Chow}, which  can be taken as a   $\cC^\infty$   version  and for  arbitrary codimension  of the proof in  \cite{FR},  which holds  under much weaker regularity assumptions. 
As we mentioned above,   the details of this   proof lead to  Corollary \ref{cor43} (see  also \cite[Thm. 4.4]{FR})  which give very useful information.   
\par
In the next  subsection \S \ref{pre-lemmas},  two preliminary lemmas on iterated Lie brackets are given. The   proof of  Theorem  \ref{Chow} is   in \S \ref{proofofchow}.
\par
\smallskip
\subsection{Two  lemmas on iterated Lie brackets}\label{pre-lemmas}
Given two smooth  (local) vector fields $X, Y$ on a manifold $\cM$, it is well known that  the commutator of their flows   satisfies
\beq \frac{1}{2} \frac{d^2}{dt^2} \left.\left(\Phi^X_{-t} \circ \Phi^Y_{t} \circ \Phi^X_{t} \circ \Phi^Y_{-t}(y)\right)\right|_{t = 0} =  - [Y, X]_y =  [X, Y]_y \qquad \text{for any}\ y\ .\eeq
This identity admits a generalisation for iterated  flow commutators which we now discuss (see also \cite{FR}). \par
\smallskip
Let     $X_1, \ldots, X_p$ be a $p$-tuple of smooth (locally defined) vector fields on an open subset $\cV \subset \cM$ and,  for any ordered $r$-tuple of integers $(i_1, \ldots, i_r)$   with $1 \leq i_1,  i_2,  \ldots ,   i_r \leq p$, let  
us denote by $X_{(i_1, \ldots, i_r)}$ the vector field defined by 
\beq \label{type-1} 
\begin{split}
&X_{(i_1)} \= X_{i_1}\qquad \text{if} \ r = 1\ ,\\
&X_{(i_1, \ldots, i_r)} \=  [X_{i_1}, [X_{i_2}, [ X_{i_3}, \ldots, [X_{i_{r- 1}}, X_{i_r}] \ldots]]]\qquad \text{if} \ r > 1\ .
\end{split}
\eeq
 For any  iterated Lie bracket  $X_{(i_1, \ldots, i_r)}$,  we    introduce the following  one-parameter family of local diffeomorphisms $G_{(i_1, \ldots, i_r) t}: \cV \subset \cM \to \cM$ which we define inductively on the integer $r$ as follows: 
\begin{itemize}
\item For $r = 1$, we set
\beq \label{G1} G_{(k)t} \= \Phi_t^{X_{k}}\ ,\eeq
\item Assuming that  $G_{(\ell_1, \dots, \ell_m)t} $ is defined for  any $1 \leq m \leq r-1$,  we set
\beq \label{G2} G_{(k, \ell_1,\ldots, \ell_{r-1})t} \=   \Phi_{-t }^{X_{k}} \circ G_{(\ell_1, \ldots, \ell_{r-1}) t} \circ  \Phi_t^{X_{k}} \circ \left(G_{(\ell_1, \ldots, \ell_{r-1}) t}\right)^{-1} \ .\eeq
\end{itemize}
By construction,  the   local diffeomorphism  $G_{(i_1, \ldots, i_r)t}$   smoothly depends on $t$. As it can be checked by  an inductive argument,  the map $G_{(i_1, \ldots, i_r)t}$  is  obtained by  composing  an appropriate  set of  flows, whose  cardinality depends on $r$ and   is equal to 
\beq \label{ndir} n_r = 2^{r-1} + \sum_{j = 1}^{r -1} 2^j = 2^r + 2^{r-1} - 2\ .\eeq 
Moreover
\begin{lem} \label{lemma4.1}
For any  $y \in \cV$  
\beq\label{theformula}
\begin{split}
& \left.  G_{(i_1, \ldots, i_r)t}(y)\right|_{t = 0} = y\ , \\
& \left. \frac{d G_{(i_1, \ldots, i_r)t}(y)}{dt}\right|_{t = 0}  =  \left. \frac{d^2 G_{(i_1, \ldots, i_r)t}(y)}{dt^2}\right|_{t = 0}  =  \ldots =  \left. \frac{d^{r-1} G_{(i_1, \ldots, i_r)t}(y)}{dt^{r-1}}\right|_{t = 0} = 0\ ,\\
&  \left. \frac{d^r G_{(i_1, \ldots, i_r)t}(y)}{dt^r}\right|_{t = 0} =  r! X_{(i_1, \ldots, i_r)}\big|_{y}\ .
\end{split}
\eeq
\end{lem}
\begin{pf}
 For $r = 1$, the equalities  follow immediately from the definitions. Assume now that the claim is true for $r = r_o -1$ and let us prove it for $r = r_o$. 
Since  the result  is of local nature, with no loss of generality we may identify $\cV $ with an open subset of $\bR^N$ and, for any  integer $v\geq 1$,   we may  consider the map  $(t, y) \mapsto G_{(\ell_1, \ldots, \ell_v)t}(y)$
as a smooth map from an open subset of  $\bR \times \bR^N$ into $\bR^N$: 
$$G_{(\ell_1, \ldots, \ell_v)}:  (-\ve, \ve) \times \cV \subset \bR \times \bR^N \longrightarrow \bR^N\ .$$ 
According to  this identification,  a vector field $X = X^i\frac{\p}{\p x^i}$ can be identified with  the $t$-independent map 
$$X:  (-\ve, \ve) \times \cV \subset \bR \times \bR^N \longrightarrow \bR^N\ ,\qquad  X \= \left(X(t, y) =  (X^1( y) , \ldots, X^N(y) \right)\ .
$$
By inductive hypothesis, 
\begin{equation} \label{prot} \frac{\p^\ell G_{(i_2, \ldots, i_{r_o})t} (y)}{\p t^\ell}  \bigg|_{t = 0}  \hskip - 10pt = 0\ , \ 1 \leq \ell \leq r_o -2 \ , \ \ \ 
 \frac{\p^{r_o -1} G_{(i_2, \ldots, i_{r_o})t} (y)}{\p t^{r_o-1}}  \bigg|_{t = 0}   \hskip - 10pt  = (r_o -1)!X_{(i_2, \ldots, i_{r_o})}(y) \ ,  \end{equation}
or, equivalently, considering the Taylor   expansion  with respect to $t$ up to order $r_o -1$ 
\beq\label{45} G_{(i_2, \ldots, i_{r_o})t} (y)  = y + X_{(i_2, \ldots, i_{r_o})}(y) t^{r_o -1} + \gr_{r_o -1}(y,t)\ ,\eeq
where  $\gr_{r_o -1}(y,t) = o(t^{r_o -1})$ for any $y$. 
  Note also that, for any given $t$, a  local inverse map 
$ \left(G_{(i_2, \ldots, i_{r_o})t}\right)^{-1} $ from an open subset of $\bR^N$ into $\bR^N$ has the form 
\beq \left(G_{(i_2, \ldots, i_{r_o})t}\right)^{-1} (y)  = y - X_{(i_2, \ldots, i_{r_o})}(y) t^{r_o -1} +  \wt \gr_{r_o -1}(y,t)\ ,\eeq
where, again,   $\wt \gr_{r_o -1} (y, t)= o(t^{r_o -1})$ for any $y$. 
 We now introduce the two-parameter family of local diffeomorphisms 
\begin{multline*}
 \wt G_{(i_1, \ldots, i_{r_o})} : (-\ve, \ve) \times  (-\ve, \ve) \times \cV\longrightarrow \bR^N\ ,\\
 \wt G_{(i_1, \ldots, i_{r_o})(t,s)} (y)\= \left(\Phi_{-t}^{X_{i_1}} \circ  G_{(i_2, \ldots, i_{r_o})s} \circ  \Phi_{t}^{X_{i_1}} \circ  \left(G_{(i_2, \ldots, i_{r_o})s}\right)^{-1} \right)(y)\ .
 \end{multline*} 
and notice that 
$$G_{(i_1, \ldots, i_{r_o})t}(y) =  \wt G_{(i_1, \ldots, i_{r_o})(t,t)} (y)\qquad \text{for any} \ y\ .$$
From  \eqref{prot}  the  derivatives in $t$ and $s$ at $(0,0)$ of the maps $\wt G_{(i_1, \ldots, i_{r_o})(t,s)}$   are such that:
\begin{itemize}[leftmargin = 20pt]
\item[(i)] For  $1 \leq j$,
\begin{multline} \frac{\p^j \wt G_{(i_1, \ldots, i_{r_o})(t,s)}( y)}{\p s^j }\bigg|_{(0,0)} = \frac{\p^j}{\p  s^j} \left( \wt G_{(i_1, \ldots, i_{r_o})(0,s)}( y)\right)\bigg|_{s = 0} = \\
= \frac{\p^j}{\p  s^j}  \left(\Id_{\bR^N} \circ G_{(i_2, \ldots, i_{r_o})s} \circ \Id_{\bR^N} \circ  \left(G_{(i_2, \ldots, i_{r_o})s}\right)^{-1}\right)(y)\bigg|_{s = 0} = \frac{\p^j}{\p  s^j} \left( \operatorname{Id}_{\bR^N} (y)\right)\bigg|_{s = 0}  = 0\ ,
\end{multline}
\begin{multline} \frac{\p^j \wt G_{(i_1, \ldots, i_{r_o})(t,s)}(y)}{\p t^j }\bigg|_{(t = 0, s= 0)} = \frac{\p^j}{\p  t^j} \left( \wt G_{(i_1, \ldots, i_{r_o})(t,0)}(y)\right)\bigg|_{t = 0} = \\
= \frac{\p^j}{\p  t^j}  \left(\Phi_{-t}^{X_{i_1}} \circ \operatorname{Id}_{\bR^N} \circ  \Phi_{t}^{X_{i_1}} \circ \operatorname{Id}_{\bR^N}\right)(y)\bigg|_{t = 0} = 
\frac{\p^j}{\p  t^j} \left( \operatorname{Id}_{\bR^N} (y)\right)\bigg|_{t = 0}  = 0\ ;
\end{multline}
\item[(ii)] For $1 \leq \ell \leq r_o -2$ and   $ j \geq 1$  
\begin{equation} \frac{\p^{j+\ell} \wt G_{(i_1, \ldots, i_{r_o})(t,s)}(y) }{\p t^j \p s^\ell }\bigg|_{(t= 0,s = 0)} =  \frac{\p^j }{\p t^j}\left(\frac{\p^\ell \wt G_{(i_1, \ldots, i_{r_o})(t,s)}(y)}{\p s^\ell }\bigg|_{(s = 0)}\right)\bigg|_{t = 0}  = 0\ ;
\end{equation}
 \item[(iii)] For  $j = 1$ and  $ \ell = r_o -1$, 
\begin{multline}  \frac{\p^{r_o} \wt G_{(i_1, \ldots, i_{r_o})(t,s)}(y) }{\p t \p s^{r_o-1} }\bigg|_{(t= 0,s = 0)} = \frac{\p}{\p t}\left(\frac{\p \wt G_{(i_1, \ldots, i_{r_o})(t,s)}(y)}{\p s^{r_o-1} }\bigg|_{(t,0)}\right)\bigg|_{t = 0} = \\
=  (r_o -1)!  \frac{\p}{\p  t}  \left((\Phi_{-t}^{X_{i_1}} )_*\left(X_{(i_2, \ldots, i_{r_o})}\big|_{\Phi_{t}^{X_{i_1}}(y)}\right) \right)\Big|_{t = 0}+\\
-(r_o -1)! \frac{\p}{\p  t}  \left((\Phi_{-t}^{X_{i_1}} \circ \operatorname{Id}_{\bR^N} \circ  \Phi_{t}^{X_{i_1}} )_* \left( X_{(i_2, \ldots, i_{r_o})}\big|_y\right)\right)\Big|_{t = 0}    =\\
=  (r_o -1)! \cL_{X_{i_1}}X_{(i_2, \ldots, i_{r_o})}\big|_{y} =    (r_o - 1)! X_{(i_1, i_2, \ldots, i_{r_o})}\big|_{y}\ ,
\end{multline}
 \end{itemize}
 where  we used the notation  $\cL_{Y} (Z) = [Y, Z]$  for  the classical Lie derivative of a vector field $Z$ along a vector field $Y$.  From (i) -- (iii), it follows that  $\wt G_{(i_1, \ldots, i_{r_o})}(t,s,y) $ has the form 
  \beq \wt G_{(i_1, \ldots, i_{r_o})}(t,s,y) = y + X_{(i_1, i_2, \ldots, i_{r_o})}(y) t s^{r_o -1} + \mathfrak R(y, t,s)\ ,\eeq
  where  $  \mathfrak R(y, t,s)$ is a smooth function of  $(y,t,s)$ which is  infinitesimal of order higher than $r_o$ with respect to $|(t,s)|$  for any given $y$. This implies that 
  \beq\label{4.13} G_{(i_1, \ldots, i_{r_o})t}(y) =  \wt G_{(i_1, \ldots, i_{r_o})}(t,t,y)  = y + X_{(i_1, i_2, \ldots, i_{r_o})}(y) t^{r_o } + \gr_{r_o}(y, t)\ ,\eeq
  and hence that \eqref{theformula} is true for $r = r_o$ as well. 
  \end{pf}
  \par
  \medskip
  Now, for any vector field $X_{(i_1, \ldots, i_r)}$,  we consider   the following associated one-parameter family of local diffeomorphisms  (defined only for  $t\geq 0$)
  \beq \label{4.13*} g_{(i_1, \ldots, i_r) t} \= G_{(i_1, \ldots, i_r)t^{\frac{1}{r}}}\ .\eeq
From  \eqref{theformula} it follows   that:  
\begin{itemize}
\item[(i)] $g_{(i_1, \ldots, i_r) t = 0}(y) = y$; 
\item[(ii)] On $t > 0$, the family is  smoothly depending on $t$  and at $t = 0$ the first  derivative  $ \frac{d g_{(i_1, \ldots, i_r)t}(y)}{dt}\Big|_{t = 0} = \lim_{t \to 0}  \frac{d g_{(i_1, \ldots, i_r)t}(y)}{dt}$ is well defined and  equal to 
\beq\left. \frac{d g_{(i_1, \ldots, i_r)t}(y)}{dt}\right|_{t = 0} = X_{(i_1, \ldots, i_r)}\big|_{y}\ ;\eeq
\item[(iii)] For   $t > 0$,  the point $g_{(i_1, \ldots, i_r)t}(y)$  and the Jacobian  
$J\left(g_{(i_1, \ldots, i_r)t}\right)\big|_{y}$ of the local diffeomorphism $g_{(i_1, \ldots, i_r)t}$ are  continuous and such that  
\beq \lim_{t = 0^+}g_{(i_1, \ldots, i_r)t}(y) = y\ ,\qquad \lim_{t \to 0^+} J\left(g_{(i_1, \ldots, i_r)t}\right)\Big|_y = \Id_{\bR^N}\ .  \eeq
\end{itemize}
With these properties,   we  are now  ready to prove the next lemma. 
 In the statement, the norms $\|\cdot\|$ of  vectors are  determined by some  choice of a Riemannian metric (it does not matter which)  on   $\cV$. 
\begin{lem} \label{lemma1} For any vector field   $ X_{(i_1, \ldots, i_r)}$, a relatively compact subset $\cV' \subset \cV$ and  a prescribed real number $\h > 0$,  there is a  $0 <  \d_o = \d_o(\h, \cV')$  and a  $1$-parameter family of local diffeomorphisms  $f^{(\d)}_t: \cV'  \to \cM$, $t \in \left(-\frac{\d}{2}, \frac{\d}{2}\right)$ associated with any $\d \in (0, \d_o)$,  smoothly  depending on $t$, which satisfies  the following conditions: 
\begin{itemize}
\item[(1)]  Each local diffeomorphism $f^{(\d)}_t: \cV' \to \cM$, $t \in \left(-\frac{\d}{2}, \frac{\d}{2}\right)$,   has  the structure
\beq \label{2.3.bis} f^{(\d)}_t = \Phi^{X_{\ell_{A}}}_{\s_{2A}(t)} \circ  \Phi^{X_{\ell_{A-1}}}_{\s_{2 A -1}(t)}  \circ \ldots \circ  \Phi^{X_{\ell_1}}_{\s_{A+1}(t)} \circ \Phi^{X_{\ell_1}}_{\s_A(t)} \circ  \Phi^{X_{\ell_2}}_{\s_{A-1}(t)}  \circ \ldots \circ  \Phi^{X_{\ell_A}}_{\s_1(t)} \ ,\eeq
for appropriate  integers $A, \ell_1, \ldots, \ell_A \geq 1$ and   real $\cC^\infty$ functions  $ \s_j(t)$; 
\item[(2)] $f^{(\d)}_{t= 0} = \operatorname{Id}_{\cV'}$; 
\item[(3)] For all $y \in \cV'$  the derivative $ \frac{d f^{(\d)}_t(y)}{dt} \Big|_{t = 0}$  is such that 
 $$\left\|  \frac{d f^{(\d)}_t(y)}{dt} \Big|_{t = 0} - X_{(i_1, \ldots, i_r)}|_y\right\| < \h\ .$$
\end{itemize}
\end{lem}
\begin{pf}
 For any sufficiently small $\d$ 
  we define  
 \beq h^{\d}_{(i_1, \ldots, i_r) t} \= g^{-1}_{(i_1, \ldots, i_r) \d}  \circ g_{(i_1, \ldots, i_r) t + \d} \ , \qquad  t \in \textstyle \big( -\frac{\d}{2}, \frac{\d}{2}\big)\ .\eeq
Note  that each    $h^{\d}_{(i_1, \ldots, i_r) t}$ satisfies the following three conditions: 
 \begin{itemize}
\item[(A)] $h^{\d}_{(i_1, \ldots, i_r) t = 0}(y) = y$; 
\item[(B)] The one-parameter family of local diffeomorphisms $t \mapsto h^{\d}_{(i_1, \ldots, i_m) t}$    is   smoothly depending on $t$; 
\item[(C)]  The first derivative with respect to $t$ at   $t = 0$  and  at the point $y$ is the vector 
\beq   \frac{d h^{\d}_{(i_1, \ldots, i_r)t}(y)}{dt}\Big|_{t = 0}  =  \bigg(g^{-1}_{i_1, \ldots, i_r) \d }\bigg)_* \left( \frac{d g_{(i_1, \ldots, i_r)t}(y)}{dt}\Big|_{t =   \d}\right)\   \eeq
and depends continuously on $\d> 0$. Moreover,  by the properties (ii) and  (iii) of the maps $g_{(i_1, \ldots, i_r) t}$ and estimating the difference
\beq
 \left. \frac{d h^{\d}_{(i_1, \ldots, i_r)t}(y)}{dt}\right|_{t = 0} - X_{(i_1, \ldots, i_r)}\big|_{y}\eeq
by means of  the integral form for the remainder  $\gr_{r}(y, t)$  in \eqref{4.13},  one can  check  the existence of a constant $C_{\cV'}$ for any relatively compact  subset $\cV' \subset \cV$   such that 
\beq  \label{limite} \left\|  \left. \frac{d h^{\d}_{(i_1, \ldots, i_r)t}(y)}{dt}\right|_{t = 0} - X_{(i_1, \ldots, i_r)}\big|_{y}\right\| < C_{\cV'}  \d^{\frac{1}{r}}\qquad \text{for any} \ y \in \cV' \ .\eeq
\end{itemize}
Furthermore, writing explicitly the local diffeomorphisms  $G_{(i_1, \ldots, i_r)t}$ in terms of the flows $\Phi^{X_{j_\ell}}_t$ and substituting such explicit expressions in the   definition of $ g_{(i_1, \ldots, i_r) t }$, one can  check that each local  diffeomorphisms $h^{\d}_{(i_1, \ldots, i_r) t} $ has  the form  \eqref{2.3.bis} with the  $\s_j(t)$  that are either smooth  functions of the form $\s_j(t) = \pm (t + \d)^{\frac{1}{r}}$ or constant functions of the form $\s_j(t) \equiv\pm  \d^{\frac{1}{r}}$. 
From this it follows that   the one-parameter family of local diffeomorphisms
$f^{(\d)}_t \= h^{\d}_{(i_1, \ldots, i_r) t}$ 
satisfies  (1) and (2). Moreover, from  \eqref{limite},  there exists   $\d =  \d(\h, \cV')> 0$   such that  $f^{(\d)}_t$ satisfies    (3) as well. 
\end{pf} 
\par
\smallskip
\subsection{Proof of  Theorem \ref{Chow}}
\label{proofofchow}
 For  any fixed  point $x_o \in \cM$,  let us denote by $\cT^{(x_o)}$   the $\cD$-path connected component of $x_o$ and  by  $\cS^{(x_o)}$   the maximal integral leaf through $x_o$  of the involutive distribution $\cD_{-\mu}$, that is the $\cD_{-\mu}$-path component of $\cM$ containing  $x_o$. Since $\cD \subset \cD_{-\mu}$, the $\cD$-path component $\cT^{(x_o)}$  is  a subset   of $\cS^{(x_o)}$. 
Being  $\cS^{(x_o)}$  (pathwise) connected, the theorem is proved if we can show that  $\cT^{(x_o)}$ is an open and closed subset of  $\cS^{(x_o)}$ in  its intrinsic topology. 
 As a matter of fact, it suffices to show that,  for any  $y \in \cS^{(x_o)} = \imath(\wt \cS)$,  the  $\cD$-path connected component   $\cT^{(y)}$ of $y$ is  an open subset of  $\cS^{(x_o)}$. This would directly imply that the  sets 
$\cT^{(x_o)}$ and  $\cS^{(x_o)} \setminus \cT^{(x_o)} = \bigcup_{y \in \cS^{(x_o)} \setminus \cT^{(x_o)}} \cT^{(y)}$
are  both open in $\cS^{(x_o)} $. 
\par
\smallskip
Let    $M \= \rank \cD_{-\mu}$. Given  $y \in \cS^{(x_o)}$,  let us consider   local generators for  the quasi-regular set $V$ of the generalised distribution $(V, \cD = \cD^V)$ nearby $y$, i.e.  a $p$-tuple of vector fields   $(X_1, \ldots, X_p)$ on a neighbourhood $\cV$ of $y$ such that any vector field $Y \in V|_{\cV}$ has the form
$$ Y  =  f^i X_i\qquad\text{for some smooth}\ \  f^i: \cV \to \bR\ .$$
As in \S \ref{pre-lemmas}, for any $1 \leq r \leq \mu$ and any ordered $r$-tuple of integers $(i_1, \ldots, i_r)$   with $1 \leq i_1 \leq i_2 \leq \ldots \leq i_{r-1} < i_r \leq p$,  we denote 
\beq \label{type-1*} 
\begin{split}
&X_{(i_1, \ldots, i_r)} \=  [X_{i_1}, [X_{i_2}, [ X_{i_3}, \ldots, [X_{i_{r- 1}}, X_{i_r}] \ldots]]]\qquad \text{if} \ r > 1\ .
\end{split}
\eeq
 By Lemma \ref{lemmino} and the hypothesis of uniform $\mu$-type,  we may consider  $M$ vector fields $Y_1, \ldots, Y_M$, each of them of the form   $Y_\ell = X_{(i_1(\ell), \ldots, i_{r_\ell}(\ell))}$ for some choice of the indices  $i_j(\ell)$,  that are  pointwise linearly independent  and  generate  the  involutive {\it regular} distribution $\cD_{-\mu}$ on the neighbourhood $\cV$ of $y$.  \par
 By  Lemma \ref{lemma1}, for any choice of   $\h > 0$ we may consider  a  relatively compact neighbourhood $\cV' \subset \cV$ of $y$, a real number $\d = \d(\h, \cV')> 0$ and a collection of  $M$ one-parameter families of local diffeomorphisms $ f^{(\d)(\ell)}_t$,  with $t$ varying in $\left(-\frac{\d}{2} , \frac{\d}{2}\right)$, each of them  defined on $\cV'$ and  satisfying (1) -- (3) for the vector fields $ Y_\ell =   X_{(i_1(\ell), \ldots, i_{r_\ell}(\ell))}$,  $1 \leq \ell \leq M$. In particular
\beq  f^{(\d)(\ell)}_{t = 0}(y) = y\ ,\qquad  \left\|\frac{d f^{(\d)(\ell)}_t(y)}{dt} \bigg|_{t = 0} -  Y_\ell|_{y}\right\|  =  \left\|\frac{d f^{(\d)(\ell)}_t(y)}{dt} \bigg|_{t = 0} -  X_{(i_1(\ell), \ldots, i_{r_\ell}(\ell))}|_{y}\right\|< \h  \\ .\eeq
  We stress the fact  that \eqref{2.3.bis} implies that, for any $1 \leq \ell\leq M$,   $t \in \left(- \frac{\d}{2}, \frac{\d}{2}\right)$ and    $z \in \cV'$, the  image $z' \= f^{(\d)(\ell)}_t(z)$   is the endpoint of a   $\cD$-path  of the form
 $$\g = \g_{1} \ast \ldots \g_{A} \ast  \g'_A  \ast \ldots \ast \g'_{1}$$
  where  the regular arcs $\g_{j}$ and $\g'_j$, $1 \leq j \leq A$,  are   integral curves of     vector fields  $X_{k_j}$  which belong to the  set of generators  $\{X_1, \ldots, X_p\}$ for the distribution $\cD$.  In particular  {\it any point $z' \= f^{(\d)(\ell)}_t(z)$ is in the same $\cD$-path connected component of $z$}. 
 \par
Now,  let $\d>0$  be a sufficiently small number  such that the composed map
\beq \label{oldmap} F^{(\d, y)}: \left(-\frac{\d}{2}, \frac{\d}{2}\right)^M \subset \bR^M \to \cM\ ,\qquad F^{(\d, y)}(s^1, \ldots, s^M) \=   f^{(\d)(1)}_{s^1} \circ \ldots \circ   f^{(\d)(M)}_{s^M}(y)\eeq
 is well defined  and of class $\cC^\infty$. By construction and the  above  remark,  we have:  
\begin{itemize}[leftmargin = 15pt]
\item[--] Any point $z = F^{(\d, y)}(s^1, \ldots, s^M)$  is the  image of the prescribed point  $y$ under a  diffeomorphism which is obtained composing maps of the form \eqref{2.3.bis} 
and hence it is connected to $y$ by a $\cD$-path;  
\item[--] In any set of coordinates,  the columns of the Jacobian  $ J\left(F^{(\d, y)}\right)\big|_{(0, \ldots, 0)}$   are given by 
the coordinate components of the vectors 
\beq \label{thevectors} \frac{d f^{(\d)(\ell)}_t(y)}{dt} \bigg|_{t = 0} = Y_\ell|_{y} + \left( \frac{d f^{(\d)(\ell)}_t(y)}{dt} \bigg|_{t = 0}  - Y_\ell\big|_{y} \right) \qquad \text{with}\ \  \left\| \frac{d f^{(\d)(\ell)}_t(y)}{dt} \bigg|_{t = 0}  - Y_\ell\big|_{y} \right \| < \h \ . \eeq
\end{itemize}
Since the  $Y_\ell|_{y}$ are   linearly independent (thus with components forming a matrix of rank $M$),  by the semicontinuity  of the rank,   for any  sufficiently small $\h$ and a corresponding  choice for   the parameter $\d = \d(\h, \cV')$,  
 the vectors $\frac{d f^{(\d)(\ell)}_t(y)}{dt} \bigg|_{t = 0} $ are linearly independent and  the  Jacobian  $ J\left(F^{(\d, y)}\right)\big|_{(0, \ldots, 0)}$  has  rank $M = \rank \cD_{-\mu}$. By the Inverse Function Theorem, this implies that, for a sufficiently small neighbourhood $ \cU \subset  \left(-\frac{\d}{2}, \frac{\d}{2}\right)^M \subset \bR^M$ 
of $0_{\bR^M}$, the image $F^{(\d, y)}(\cU)$ is an $M$-dimensional  submanifold of $\cM$, which contains $y = F^{(\d, y)}(0)$ and is entirely included in $\cT^{(y)} \subset \cS^{(y)} = \cS^{(x_o)}$.   If $\cU$ is sufficiently small,    $F^{(\d, y)}(\cU)$ is  open in one of the  sets $\cW = \imath(\wt \cW) \subset \cS^{(x_o)}$ of the intrinsic topology of $\cS^{(x_o)} = \imath(\wt \cS)$
which are embedded submanifolds of $\cM$.  By the arbitrariness of    $y$ and the remarks at the beginning of the proof,  it follows that $\cT^{(x_o)}$ is open and closed  in   the intrinsic topology of $\cS^{(x_o)} $.\hfill\square
\par
\smallskip
\section{A  corollary of the  proof} \label{corchow}  In the next corollary, we give the  information,  which  can be derived from the above   proof of Theorem \ref{Chow} and which we advertised at the beginning of \S \ref{appendix}.  But before  stating  it,   a     preliminary clarification is in need.  As we pointed out in the previous section, for any fixed point $y \in \cS^{(x_o)}$ and any choice of a  sufficiently small $\d$ (i.e. of a sufficiently small  domain), the map $F^{(\d, y)}$
takes values in a relatively compact subset  $\cW = \imath(\wt \cW)$ of $\cS^{(x_o)}$, which  contains $y$ and  is   an $M$-dimensional  embedded submanifold of $\cM$. Such $\cW$ can be assumed to be contained in a  prescribed relatively compact  neighbourhood $\cV \subset \cM$ of $y$. 
With  no loss of generality,  we may   assume that the neighbourhood  $\cV$ admits   a set of  coordinates for the $N$-dimensional manifold $\cM$  with the properties that: (a)  $y$ corresponds to  the origin $y =  0_{\bR^N}$; (b)   $\cW$ is the $M$-dimensional submanifold of $\cV$  obtained by imposing each of the last $N - M$ coordinates equal to  $0$. As a consequence,   the  first  $M$ coordinates of this chart  represent also a set of  coordinates for the submanifold  $\cW$.  From now on, {\it we  use such coordinates  to identify $\cW$ and $\cV$ as open subsets of  $ \bR^M$ and $\bR^N$, respectively}.  In this way the map  $F^{(\d, y)}$ can be  taken as a smooth map between two open subsets of $\bR^M$, which maps $0$ into itself, and we may use the standard Euclidean metric of $\bR^M$ 
to determine  distances between points and  norms  of vectors.\par 
\smallskip
We now claim that,  after fixing a point $y_o \in \cV$, there exists a sufficiently small compact neighbourhood $\cV_o \subset \cV$ of $y_o$ and a positive real number  $\d_{\max}  = \d_{\max} (y_o,\cV)> 0$, depending on $y_o$ and 
$\cV$,  such that {\it for  any $\d \in (0, \d_{\max})$,  each corresponding map  $F^{(\d, y)}$, $y \in \cV_o$,  is well defined on the ipercube $\big[- \frac{\d}{2}, \frac{\d}{2}\big]^M$}. This can be checked as follows. 
We recall that each   map  $F^{(\d, y)}$  is  obtained  by applying to $y$  a composition of    the  local diffeomorphisms  $f^{(\d)(i)}_{s}$ with an appropriate sequence of flow parameters. The latter are in turn  appropriate finite compositions of the  flows of the vector fields $X_i$, so that   $F^{(\d, y)}$ itself is determined by applying to $y$ a  finite composition  of flows. 
We now observe that  there exists   at least one  interval of definition $  I_{\text{max}} = [-\ve_{\text{max}} , \ve_{\text{max}} ]\subset \bR$  which is common to all   integral curves  $t \mapsto \g(t)$ of  the vector fields $X_i$  that start from  points   of  the relatively compact open set  $\cV$.  There is also a possibly smaller relatively compact neighbourhood  $\cV'$ of $y_o$  and a  smaller interval $I'_{\text{max}} = [-\ve'_{\text{max}} , \ve'_{\text{max}} ] \subset I_{\text{max}}$ 
 with the property that   all of the above  integral curves starting from points in $\cV' \subset \cV$ have their final points in $\cV$. By construction, the  composition  $\Phi^{X_i}_s \circ \Phi^{X_j}_{\wt s}(y)$ is  well defined for any $y \in \cV'$ and $s, \wt s \in [-\ve'_{\text{max}} , \ve'_{\text{max}} ]$. Iterating this argument,  after a finite number of reductions $\cV \supset \cV' \supset  \cV'' \supset \ldots$  and $ I_{\text{max}} \supset I'_{\text{max}}\supset I''_{\text{max}}\supset \ldots$  of the open set  $\cV$ and of the interval $I_{\text{max}}$,  one ends up with a neighbourhood $\cV_o$ of $y_o$
and a real number $\wc \ve_{\text{max}}$ such the composition of the flows which appears in the definitions of the  maps  $F^{(\d, y)}$, $y \in \cV_o$,  is well defined for the  flow parameters  running in the interval  $[-\wc \ve_{\text{max}} , \wc \ve_{\text{max}} ]$. We finally recall that, in the definition of  $F^{(\d, y)}$ the flow  parameters are either equal to a constant (which depend continuously on $\d$) or are given by  continuous functions 
of the arguments $(s^i) = (s^1, \ldots, s^M)$ of $F^{(\d, y)}$. This implies the existence of a real number
$\d_{\max}   = \d_{\max} (\cV, y_o)> 0$ such that if $\d \in (0, \d_{\max})$,  all flow parameters occurring in the construction of $F^{(\d, y)}$ take  values in  $[-\wc \ve_{\text{max}} , \wc \ve_{\text{max}} ]$ whenever  the arguments $(s^i)$ of   $F^{(\d, y)}$ are in  the hypercube $\big[- \frac{\d}{2}, \frac{\d}{2}\big]^M$.  This implies that $\d_{\text{max}}$ and $\cV_o$ have the required property.  With no loss of generality, we may also  assume    $\d_{\max} < 1$. 
 \par
In the next statement and throughout its proof, for any relatively compact set $\cV_o \subset \cM$ and any function $f: \cV_o \to \bR$ which is $\cC^k$ up to the boundary,  we denote 
$$\vertiii{ f }_{\cC^k(\cV_o)} \= \max\{ 1 \ , \|f\|_{\cC^k(\cV_o)}\}\ .$$
 \begin{cor} \label{cor43} Let $ \{Y_1, \ldots, Y_M\}$  be a fixed choice of generators for the distribution 
 $\cD_{-\mu}|_\cV$ of the form   $Y_\ell = X_{(i_1(\ell), \ldots, i_{r_\ell}(\ell))}$, with $(X_i)$ generators for $\cD$, 
 and denote by  $\cV_o$ and  $0 < \d_{\max}  = \d_{\max}(\cV) < 1$  a relatively compact neighbourhood of a point $y_o \in \cV \subset  \cM$ and a corresponding  real number  such that the  property described above   holds. 
 \par
Assume that  $\bigg|\det \left(Y_i^j(y)\right)_{1 \leq i, \ell \leq M}\bigg|$ is bounded away from $0$  in $ \cV_o$,  Then there are  constants $K, K' > 0$,  depending just on  $M = \dim \cW$,  $N = \dim \cM$ and $\mu$,  such that for any $y \in \cV_o$ and  any 
 $C_0, C_1> 0$   satisfying  
 \beq \label{thebounds} C_0 \leq \inf_{y \in \cV_o} \left|\det \left(Y_i^j(y)\right)_{1 \leq i, \ell \leq M}\right| \ ,\qquad  \sum_{j,\ell} \vertiii{X^j_\ell}_{\cC^{\mu+3}(\cV_o)}\leq C_1\ ,\eeq 
   the  image of the map  $F^{(\d_o, y)}: \left(-\frac{\d_o}{2} , \frac{\d_o}{2}\right)^M \subset \bR^M \to \cM$   contains the ball $B_{\mathbf r_o}(y)$,  centred at $ y \equiv 0_{\bR^m}$ and of radius  
 \beq \label{radietto} \mathbf r_o \=   K' \d_o \frac{\min\{C_0, C_0^2\}}{C_1^{{\mathbf N}(\mu, M)}} \ ,\eeq
 where \beq \label{thedo} \d_o \= \min\left\{K \frac{C_0^\mu} { C_1^{\mu{\mathbf N}(\mu, M)}}, \d_{\max} \right\}\  , \ \text{with}\ {\mathbf N}(\mu, M){\=}{{6 (2^\mu + 2^{\mu -1} - 2) (\mu + 3)M(2 M +1)}}.\eeq
 \end{cor}
  \begin{rem}   Corollary \ref{cor43}   shows  that,  if there are upper and lower  bounds as  in \eqref{thebounds}   for any $ y_o \in  \cM$, then for any point $y_o$ of the manifold there exist a neighbourhood $\cV_o$ and a real number $\d_o > 0$  such that  all maps $F^{(\d_o, y)}$,  $y \in {\cV_o}$, are surjective onto   balls, each of them with   a  radius ${\mathbf r}_o$ that  is  independent  of  $y$.  This is the property that motivates our interest. In fact,   under these hypotheses,  if  $y_k $  is a  sequence   converging   to a point   $y_o$, then  for any sufficiently large $k$   the image of the  corresponding  map $F^{(\d_o, y_k)}$  contains a whole neighbourhood of $y_o$.  Similarly,  given a sequence of uniform $\mu$-type distributions $\cD^{(k)}$ on $\cM$,  admitting a set of  generators   $\{Y^{(k)}_1, \ldots,  Y^{(k)}_M\}$ that converge  (in an appropriate norm)  to a set of  generators $\{Y^{(\infty)}_1, \ldots,  Y^{(\infty)}_M\}$  for a limit distribution $\cD^{(\infty)}$,  a pair of uniform bounds  as  in \eqref{thebounds} for the distributions  $\cD^{(k)}$ imply  the  existence of  balls of  appropriate  radii in the images of the  maps $F^{(\d_o, y)} $  of  the limit distribution $\cD^{(\infty)}$.  
  \end{rem}
 \begin{pf} Let $\h$ be the real number $\h = \frac{1}{2}  \det \left(Y_i^j|_{y = 0}\right)_{1 \leq i, \ell \leq M} $. Via a reordering of the generators, we may assume that $\h> 0$. Let us  use  such a value  in the   construction of   the map $F^{(\d, y)}$ given   in  \S \ref{proofofchow}. We recall that  for any $0 < \d \leq \d_{\max}$ and any $y \in \cV_o$, the  corresponding map $F^{(\d, y)}$ is well defined on $[- \frac{\d}{2}, \frac{\d}{2}]^M \subset \bR^M$ and is determined as  an appropriate composition of  maps 
  \beq \label{exp1}   f^{( \d)(\ell)}_{ s}(z) \=  g_{I_\ell\, \d}^{-1} \circ g_{ I_\ell\,s + \d}(z) \ ,\qquad 1 \leq \ell \leq M\ ,\eeq
 where:  
 \begin{itemize}[leftmargin = 10pt]
\item $I_\ell \= (i_1(\ell), \ldots, i_{r_\ell}(\ell))$ denotes the tuple of integers  such that   $ Y_\ell =  X_{I_\ell} =  X_{(i_1(\ell), \ldots, i_{r_\ell}(\ell))}$,  
\item $ g_{ I_\ell\,t} =  G_{I_\ell t^{\frac{1}{r_\ell}}}$  where  $G_{I_\ell s}$ is the composition of  flows of the vector fields $X_j$,  given  by \eqref{G1} and \eqref{G2}.
\end{itemize}
Moreover, we may take 
 $\d = \d(\h)$  sufficiently small  (with no loss of generality,  we may also assume $\d\leq \h$) such that 
\beq \label{conditi} \left\| \left. \frac{d f^{( \d)(\ell)}_{ s}(y)}{d s} \right|_y - Y_\ell\big|_{y} \right\| = \left\|\left. \frac{d  g_{I_\ell\, \d}^{-1}\circ g_{ I_\ell\, s + \d}(y)}{ds}\right|_{s = 0} - Y_\ell\big|_{y} \right\| < \h =  \frac{1}{2} \det \left(Y_i^j|_{y = 0}\right)_{1 \leq i, \ell \leq M} \ .\eeq
We want to determine an estimate for such a $\d = \d(\h)$. This can be  derived from   
a common upper bound for the second derivatives with respect to $s$ of the  maps \eqref{exp1}  at  the  points $z \in {\cV_o}$. This is what  we are now going to compute. 
By   \eqref{45}  the map 
$g_{I_\ell\, s + \d}(z) = G_{ I_\ell\, (s + \d)^{\frac{1}{r}}}(z)$   admits the expansion 
%
%
%
\beq \label{425} g_{ I_\ell\, s + \d}(z) = z + X_{I_\ell}(z) (s + \d) + {\mathfrak r}_{r_\ell}(z, (s  + \d)^{\frac{1}{r_\ell}})\ ,
\eeq
where the map $(z, s) \mapsto {\mathfrak r}_{r_\ell}(z, (s  + \d)^{\frac{1}{r_\ell}})$ can be explicitly determined from the remainder  $\gr_{r_\ell}(z,t)$   in \eqref{4.13}, written in integral form. \par
We now observe  that, for any $z \in {\cV_o}$, being each map  $t \mapsto \Phi^{X_\ell}_{ t}(z)$ an  integral curve  of the vector field $X_{\ell} $, each of its  partial  derivatives  at $z$ with respect to  $t$ and/or the coordinates $z^i$ 
is    equal   either to  a polynomial function (of  order less than or equal to the order of the considered derivative) of the coordinate components  of the vector field $X_{\ell} $ and their
 partial derivatives,  or to  the value of a solution of an ordinary differential equation with independent variable $t$  and coefficient given by the same partial derivatives of the components of $X_\ell$. In all cases, upper bounds for such partial derivatives of order $k$  can be  determined by  powers of upper bounds for the partial derivatives of the components of the $X_\ell$ up to the same order $k$. Assuming that the latter upper bounds are greater than $1$,  up to a multiplicative constant the power $k$ of these upper bounds  give a bound for  all  partial derivatives up to $k$. On the other hand, each map $(t, z) \mapsto G_{I_\ell t}(z)$
is a composition of     flows    $(t', z') \mapsto  \Phi^{X_\ell}_{ t'}(z')$, whose number is   less than or equal to $n_\mu \= 2^\mu + 2^{\mu-1} - 2$ (see \eqref{ndir}). 
Combining these two observations it follows that all derivatives of $k$-th order, $k \geq 1$,  of the maps $(t, z) \mapsto G_{I_\ell t}(z)$ with respect to $t$ and  with respect to the coordinates $z^i$,  are completely determined by  the  coordinate components  of the vector fields $X_{\ell} $ and of   their partial derivatives up to order $ k$. More precisely, upper bounds for such $k$-th order derivatives are determined (up to constants) by powers of order less than or equal to $k n_\mu + 1 - k \leq n_\mu k$ of upper bounds greater than $1$ for the derivatives  up to order $k$ of the coordinate components of the vector fields $X_i$.\par
\par
\smallskip
Writing  the rest $ \gr_{r_\ell}$ in  \eqref{4.13} in its integral form, 
one gets 
\beq \label{4545}  \gr_{r_\ell}(z,(s  + \d)^{\frac{1}{r_\ell}}) \= \frac{1}{r_\ell!} \int_0^{(s  + \d)^{\frac{1}{r_\ell}}} ((s  + \d)^{\frac{1}{r_\ell}} -\s)^{r_\ell} \frac{\p^{r_\ell + 1 } G_{I_\ell\,\s}(z)}{\p \s^{r_\ell +1}}  d \s\ .\eeq
It follows that for any $z \in {\cV_o}$ and  $s \in \left(- \frac{\d}{2}, \frac{\d}{2}\right)$,  
\begin{multline} \label{4.28} \frac{d^2 g_{ I_\ell\, s +\d}(z) }{ds^2} = \frac{1}{r_\ell!} \int_0^{(s  + \d)^{\frac{1}{r_\ell}}} \frac{d^2((s  + \d)^{\frac{1}{r_\ell}} -\s)^{r_\ell}}{d s^2} \frac{\p^{r_\ell + 1 } G_{I_\ell\,\s}(z)}{\p \s^{r_\ell +1}}  d \s = \\
=   \frac{1}{r_\ell!}  \frac{r_\ell-1}{r_\ell} \int_0^{(s  + \d)^{\frac{1}{r_\ell}}}\hskip - 0.2cm \s \left( (s + \d)^{\frac{1}{r_\ell}}- \s\right)^{r_\ell -2} (s + \d)^{\frac{1}{r_\ell} -2} \frac{\p^{r_\ell + 1 } G_{I_\ell\,\s}(z)}{\p \s^{r_\ell +1}}  d \s \leq\\
\leq
\frac{1}{r_\ell!}  \frac{r_\ell-1}{r_\ell}  \sup_{z \in \cV_o} \left| \frac{\p^{r_\ell + 1 }G_{I_\ell\,\s}(z)}{\p \s^{r_\ell +1}}  \right| \left( \frac{3\d}{2} \right)^{\frac{r_\ell -2}{r_\ell}}\left( \frac{\d}{2}\right)^{\frac{1}{r_\ell} - 2} \int_0^{(s  + \d)^{\frac{1}{r_\ell}}}\hskip - 0.2cm \s d \s \ .\end{multline}
Since    $r_\ell \leq \mu$,  this implies  that  there exists  a constant $K_1 > 0$ (depending only on $M$, $\mu$ and $N$) such that 
\beq  \label{4.29} \left\|\frac{d^2 g_{ I_\ell\, s + \d}(z) }{ds^2}\right\| <  K_1\left( \sum_{i,\ell} \vertiii{ X^j_\ell}_{\cC^{\mu +1}({\cV_o})} \right)^{n_\mu (\mu +1)}\frac{1}{\d^{1 - \frac{1}{\mu}}}  \ .
\eeq
We  observe that, being each map $g_{ I_\ell\,t}$ a composition of flows, also their inverse maps, $g^{-1}_{ I_\ell\,t}$  are composition of flows (one is obtained from the other by reversing the 
order of the  considered flows and changing all parameters  of  the flows into  their opposites).  Then, using once again the expansion \eqref{425}, the components $J\left(g^{-1}_{I_\ell\,\d}(z)\right)^i_j$  and  $J\left(g_{I_\ell\,\d}(z)\right)^i_j$ of the  Jacobian matrices  of the maps $z \to g^{-1}_{I_\ell\,\d}(z)$   and $z \to g_{I_\ell\,\d}(z)$, respectively,  admit a common upper bound   of the form 
\beq \label{4.29bis} J \left(  g^{-1}_{I_\ell\,\d}(z)\right){}^i_j\ ,\  J \left(  g_{I_\ell\,\d}(z)\right)^i_j \leq  \d^i_j +  K_2 \left(\sum_{i,\ell} \vertiii{ X^j_\ell}_{\cC^{\mu + 2}({\cV_o})}\right)^{n_\mu( \mu + 2)} \ \eeq
for some  constant $K_2 > 0$ (which, as before,  depends just on $M$, $\mu$ and $N$). 
\par
\smallskip
Combining \eqref{4.29} with \eqref{4.29bis} and with the definition of the map $f^{( \d)(\ell)}_{ s}(z)$, we infer   the existence of a constant $K_3 > 0$  (depending just on   $M$, $\mu$ and $N$),   which  allows to bound  the  second derivatives of the $f^{( \d)(\ell)}_{ s}(z)$ with  respect to $s$ or with respect to  $\d$ and $s$ at any   $s_o\in \left(- \frac{\d}{2}, \frac{\d}{2}\right)$ and  $z \in {\cV_o}$   by
\beq \label{430}  \left\|\left.\frac{d^2  f^{( \d)(\ell)}_{ s}(z) }{d\d d s}\right|_{s = s_o} \right\|\ ,\  \left\|\left.\frac{d^2  f^{( \d)(\ell)}_{ s}(z) }{ds^2}\right|_{s = s_o} \right\| <   K_3 \frac{\left( \sum_{i,\ell} \vertiii{ X^j_\ell}_{\cC^{\mu +3}({\cV_o})} \right)^{3 n_\mu(\mu + 3)}}{\d^{1 - \frac{1}{\mu}}} \ .
\eeq
 From this  we get the following.  For any $\d \leq \d_{\max}$:  
\begin{itemize}[leftmargin = 20pt]
\item   There exists a constant $K_4 > 0$,   depending only  on $M$, $\mu$,  $N$,  such that 
\begin{multline*}
\left\| \left. \frac{d f^{( \d)(\ell)}_{ s}(y)}{d s} \right|_y - Y_\ell\big|_{y} \right\| \leq  K_4   \frac{\left( \sum_{i,\ell} \vertiii{ X^j_\ell}_{\cC^{\mu +3}({\cV_o})} \right)^{3 n_\mu(\mu + 3)}}{\d^{1 - \frac{1}{\mu}}}   \d <   \\
< K_4  \left( \sum_{i,\ell} \vertiii{ X^j_\ell}_{\cC^{\mu +3}({\cV_o})} \right)^{3 n_\mu(\mu + 3)}\d^{\frac{1}{\mu}}\ ;\end{multline*}
This in turn implies that the inequality \eqref{conditi} holds for any $\d$ satisfying 
\beq \label{ulla} \d \leq  \left(\frac{1}{2 K_4}\right)^\mu  \left( \frac{\det \left(Y_i^j|_{y = 0}\right)_{1 \leq i, \ell \leq M} }{ \left( \sum_{i,\ell} \vertiii{ X^j_\ell}_{\cC^{\mu +3}({\cV_o})} \right)^{3 n_\mu(\mu + 3)}} \right)^\mu\ .\eeq
\item  There exists    a constant $K_5  > 0$ such that for any two points $\wc s  = (s^1, \ldots, s^M)$,  $\wc s' = (s'{}^1, \ldots, s'{}^M) $ on which the map $F^{(\d, y)}$ is defined, the following 
estimate holds (we  stress  that   the exponent $6 n_\mu( \mu + 3) M$ appearing in the  estimate might not be optimal; we also remind that, since  $\d \leq \d_{\max} <  1$, then  
$\d < \d^{1 - \frac{1}{\mu} } $, a property which is used to derive the second inequality) 
\begin{multline} \label{4.39} \left\| JF^{(\d, y)}\big|_{\wc s} - JF^{(\d, y)}\big|_{\wc s'}\right\| \leq K_5  \frac{\left( \sum_{i,\ell} \vertiii{ X^j_\ell}_{\cC^{\mu +3}({\cV_o})}\right)^{6 n_\mu(\mu + 3) M}}{\d^{1 - \frac{1}{\mu}}} \left| \wc s - \wc s'\right| < \\
<  K_5  \left( \sum_{i,\ell} \vertiii{ X^j_\ell}_{\cC^{\mu +3}({\cV_o})}\right)^{6 n_\mu( \mu + 3) M} \sqrt{M}   \ .\end{multline}
\end{itemize}
Let us  use this estimate to determine an (independent of $\d$ and $y$) upper bound for the component of the  matrix valued map  $\d \to  (JF^{(\d, y)}|_{\wc s = 0})^{-1}$, $\d > 0$. 
For this, we  
 recall that  the components of  $\d \to  (JF^{(\d, y)}|_{\wc s = 0})^{-1}$ are rational functions of the components of  $\d \to  (JF^{(\d, y)}|_{\wc s = 0})$, whose denominator is  $\det\left(JF^{(\d, y)}|_0\right)$ and the numerator is a polynomial of degree $M -1$.  Hence, if we denote by $D> 1$  an upper bound for the components of   $\d \to  (JF^{(\d, y)}|_{\wc s = 0})$, then,  {\it up to a multiplicative constant}, 
 the components of the map  $\d \to  (JF^{(\d, y)}|_{\wc s = 0})^{-1}$  are bounded above by    $\frac{D^{M-1}}{C_0}$.
Using the same circle of ideas,    which  led to  the estimates \eqref{430} and \eqref{4.39},   it is possible to  check  that   one can assume  
\beq D = \wt K  \left( \sum_{i,\ell} \vertiii{ X^j_\ell}_{\cC^{\mu +3}({\cV_o})} \right)^{6 n_\mu (\mu + 3) M} \eeq
 for some   constant $\wt K$ depending only on $N, M$ and $\mu$.  
Thus, setting $K_6 \= \wt K^{M-1}$, we get  that   the components of the map  $\d \to  (JF^{(\d, y)}|_{\wc s = 0})^{-1}$    are bounded above by 
\beq    \| (JF^{(\d, y)}|_{\wc s = 0})^{-1}\|  < \frac{ K_6}{C_0}  \left( \sum_{i,\ell} \vertiii{ X^j_\ell}_{\cC^{\mu +3}({\cV_o})} \right)^{6 n_\mu (\mu + 3) M^2} \hskip -1 cm <  K_6\frac{ C^{6 n_\mu (\mu + 3)M^2}_1}{C_0}\ .\eeq
We now recall that, from  standard facts on the Inverse Function Theorem (see e.g. \cite[Thm. 2.9.4]{HH}), if ${\bf r}>0$ is a real number such that 
\begin{itemize}[leftmargin=20pt]
\item[(1)] the ball of radius $2 {\mathbf r} \| (JF^{(\d, y)}|_{\wc s = 0})^{-1}\|$ and centre $\wc s  = 0$ is contained in the domain of $F^{(\d, y)}$ (i.e.  if $2 {\mathbf r} \| (JF^{(\d, y)}|_{\wc s = 0})^{-1}\| \leq \frac{\d}{2}$); 
\item[(2)]  the Lipschitz ratio 
$\frac{\left\| JF^{(\d, y)}\big|_{\wc s} - JF^{(\d, y)}\big|_{\wc s'}\right\|} {|\wc s - \wc s'|} $ is less than or equal to $\frac{1}{2 {\mathbf r}  \| (JF^{(\d, y)}|_{\wc s = 0})^{-1}\|^2} $
\end{itemize}
then the image of $F^{(\d, y)}$ contains a ball of radius ${\mathbf r}$ centred at $y$. 
From the  estimate \eqref{4.39} and denoting  ${\mathbf N}(\mu, M) \= {6 n_\mu (\mu + 3)M(2 M +1)} $, we see that  (1) and (2) are necessarily  satisfied  if  ${\mathbf r}$  is such that 
\beq \label{ineq}   {\mathbf r}\leq   \frac{1 } {  4 K_6  } \frac{\d \,C_0}{C_1^{{\mathbf N}(\mu, M)}} \ \left( \ <   \frac{1 } {  4 K_6  } \frac{\d\,C_0}{C_1^{6 n_\mu(\mu + 3) M^2}} \right) \ ,\qquad   \mathbf r   \leq  \frac{1}{ 2 K_5 K^2_6 }\frac{\d^{ \frac{\mu -1}{\mu}}C^2_0}{C_1^{ {\mathbf N}(\mu, M)}} \ .\eeq
If we set (we recall that $\d_{\max} < 1$)
\beq 
 \d_o =  \min\left\{  \d_{\max} ,  \ \ \frac{1}{(2 K_4)^\mu}\frac{ C_0^\mu }{ C_1^{\mu {\mathbf N}(\mu, M)}} \right\}\ ,\qquad
\mathbf r_o =  \min\left\{ \frac{1}{ 4 K_6},  \frac{1}{ 2 K_5 K^2_6 }\right\} \d_o \frac{\min\{C_0, C_0^2\}}{C_1^{{\mathbf N}(\mu, M)}}\ ,
\eeq
then $\d_o$ is   less than or equal to $\d_{\max}$ and  verifies \eqref{ulla},  both  (1) and (2) are satisfied by   $(\mathbf r = \mathbf r_o, \d = \d_o)$ and the ball  $B_{\mathbf r_o}(0)$ 
  is   contained in   $F^{(\d_o, y)}\big(\!\!\left(-\frac{\d_o}{2} , \frac{\d_o}{2}\right)^M\!\!\big)$. Denoting 
  $K{\=}\frac{1}{(2 K_4)^\mu} $,  $K'{\=}\min\left\{ \frac{1}{ 4 K_6},  \frac{1}{ 2 K_5 K^2_6 }\right\}$, 
  the claim follows.  
 \end{pf}
 \par
 \bigskip
\begin{example}   A hint of  the possible applications of  Corollary \ref{cor43}  might be obtained  from  the following simple  example. 
 Consider the   $1$-parameter family of  pairs of   vector fields $(X_1^{(\l)}, X_2^{(\l)})$, $\l \in \bR$,  defined   on  $\bR^3$  by 
\beq \label{ecco11}  X^{(\l)}_1 \= \frac{\p}{\p x^1}\ ,\qquad X^{(\l)}_2 \= \frac{\p}{\p x^2} + \bigg(x^1 +  \l  (x^1)^5 \sqrt{|x^1|} \r(x^1, x^2, x^3) \bigg) \frac{\p}{\p x^3}\eeq
where we denote by   $\r: \bR^3 \to [0, 1]$  a  $\cC^\infty$  bump function  which is identically equal to $1$ on the unit closed ball $\overline{B_1(0)}$ and with  support in the ball $ B_2(0)$ of radius $2$. 
For any $\l$ let  us denote by $(V^{(\l)}, \cD^{(\l)})$   the  pair given by the family $V^{(\l)}$ of the local vector fields  of $\bR^3$  of the form $Y \= f_1 X_1^{(\lambda)} + f_2 X_2^{(\lambda)} $,  determined by  $\cC^\infty$ functions   $f_1$, $f_2$,  and the corresponding family of  tangent subspaces,  defined  at each point by 
$$\cD^{(\l)}_{x} = \big\langle X^{(\l)} _1|_x, X^{(\l)}_2|_x \big\rangle\ .$$
Note that   $\cD_o \= \cD^{(\l = 0)}$  is the standard  {\it contact structure} of $\bR^3$.  It  is  a  regular $\cC^\infty$ distribution  of uniform $\mu$-type with $\mu = 2$, 
and its  associated distribution $\cD^{(\text{Lie})}_o$ coincides with $T \bR^3$, being    generated by the  vector fields 
$$Y_1 \= X_1^{(\l= 0)} = \frac{\p}{\p x^1}\ ,\quad Y_2 \= X^{(\l= 0)}_2 \= \frac{\p}{\p x^2}  + x^1 \frac{\p}{\p x^3} \  ,\quad  Y_3 = [X^{(\l= 0)}_1,  X^{(\l= 0)}_2] = \frac{\p}{\p x^3} \ .$$
As  immediate consequence of Theorem \ref{Chow}, any two points $x, x' \in \bR^3$ can be joined by a $\cD_o$-path, as  is well known.\par
\smallskip
On the other hand,   none of the pairs $(V^{(\l)}, \cD^{(\l)})$,  $\l \neq 0$,  fits into the class considered in Definition \ref{quasi-regular-def}, because 
 one of the two generators   is   not  even    $\cC^6$.  Nonetheless,  it still   makes sense to consider  $\cD^{(\l)}$-paths (=  piecewise regular curves  that  are tangent to  the spaces of 
$\cD^{(\l)}$  at  all points) and to ask  which  pairs of  points   can be joined    by   $\cD^{(\l)}$-paths. 
Despite of the fact that Theorem \ref{Chow} cannot be immediately  used for  the pairs $(V^{(\l)}, \cD^{(\l)})$,  
using  Corollary \ref{cor43} one can   prove    that {\it  any   pair  $x, x' \in \bR^3$  can  be   joined by a $\cD^{(\l)}$-path for any   $\l$  in a sufficiently small neighbourhood of $0$}.  
\par
To see this,   consider the triples $(Y_i^{(\l)})_{i = 1,2,3}$ of vector fields (depending on $\l$) 
\beq Y^{(\l)}_1 \= X_1^{(\l)} \ ,\quad Y^{(\l)}_2 \= X^{(\l)}_2\  ,\quad  Y^{(\l)}_3 \= [X^{(\l)}_1,  X^{(\l)}_2]\ . \eeq
Using just the definition, one can directly see that, for any  $\ve > 0$, there exists a    $\d > 0$     such  that the coordinate components $X_i^{(\l)j}$, $Y_{\ell}^{(\l)j}$, of the vector  fields $X_i^{(\l)}$, $Y_\ell^{(\l)}$
satisfy 
\beq \max_{\smallmatrix i = 1,2 \\ j = 1,2,3
\endsmallmatrix}\|X_i^{(\l)j} - X_i^{(\l = 0)j}\|_{\cC^5(\overline B_3(0))}, \quad  \max_{\smallmatrix \ell = 1,2, 3\\ j = 1,2,3
\endsmallmatrix}   \|Y_{\ell}^{(\l)j} - Y_{\ell}^{(\l = 0)j}\|_{\cC^0(\overline B_3(0))} < \ve \eeq
for any $\l \in (-\d, \d)$.
Then, using the   smoothness of the components of the  vector fields  $X_i^{(\l = 0)}$, $Y_\ell^{(\l = 0)}$, we have that if $\ve $ is sufficiently small there are   constants 
$C_0$, $C_1$ such that  the estimates  \eqref{thebounds} are satisfied  over the ball $\cV_o = \overline{B_3(0)} \subset \bR^3$ 
by the coordinate components  of the  vector fields $X_i^{(\l)}$, $Y_{\l}^{(\ell)}$  for any $\l \in (-\d, \d)$.  \par
Now, take  a sequence of $\cC^\infty$  functions $f_k: [-2, 2]\subset \bR \to \bR$,  converging  in   $\cC^5$-norm of $[-2, 2]$ to the function 
$f(s) = s^5 \sqrt{|s|}$  for $k \to \infty$ (a sequence  of this kind can be  determined  using  Bernstein polynomials -- see e.g. \cite{Ph}).  Then, for each $\l \in (- \d, \d)$ and $k \in \bN$, consider the vector fields $X^{(\l|k)}_i$, $Y^{(\l|k)}_\ell$ 
which can be obtained from the vector fields  $X^{(\l)}_i$, $Y^{(\l)}_\ell$ by simply replacing  each occurrence of the  expression $f(x^1) = (x^1)^5 \sqrt{|x^1|}$ by   $f_k(x^1)$. The new vector fields $X^{(\l|k)}_i$ (resp. $Y^{(\l|k)}_\ell)$  are  $\cC^\infty$ and their coordinate components  converge  to the components of  the vector fields $X^{(\l)}_i$ (resp. $Y^{(\l)}_\ell)$  in  the  $\cC^5$- (resp.  $\cC^4$-) norm on  the compact set $\overline{B_3(0)} \subset \bR^3$. Thus, after possibly slightly adjusting the constants $C_0$ and $C_1$, for a fixed $\l    \in (- \d, \d)$, the estimates  \eqref{thebounds} hold  on $\cV_o = \overline{B_3(0)}$  for  $X^{(\l|k)}_i$, $Y^{(\l|k)}_\ell$, $k \geq k_o$,   with  $k_o$ sufficiently large. \par
By Corollary \ref{cor43}, for   any $y \in  \overline{B_2(0)}$ and for a possibly larger  $k_o$,  for each $k \geq k_o$ there is a  map $F^{(k)(\d_o, y)}: \left(-\frac{\d_o}{2} , \frac{\d_o}{2}\right)^3 \subset \bR^3 \longrightarrow  \bR^3$
with the following three properties: 
\begin{itemize}[leftmargin = 20pt]
\item[(i)] Its domain is an hypercube  $\left(-\frac{\d_o}{2} , \frac{\d_o}{2}\right)^3 \subset \bR^3$, which is independent of   $k$ and   $y$; 
\item[(ii)] Its  image    contains   a  ball $B_{\mathbf r_o}(y)$, whose  radius  ${\mathbf r_o}$ is  independent of $k$ and  $y$; 
\item[(iii)] The points   $F^{(k)(\d_o, y)}(s)$ are the images of $y$ under an appropriate composition   (defined  in \eqref{oldmap} and  \eqref{2.3.bis})  of flows of  the vector fields $X_i^{(\l|k)}$, which  can be  assumed to be determined by a  construction that    is independent of  $k$ and $y$.  
\end{itemize}
This  imply that, for  each $y \in \overline{B_2(0)}$ and  for $k \to \infty$,  the  smooth maps   $F^{(k)(\d_o, y)}$ converge  in $\cC^5$-norm to a map  $F^{(\d_o, y)}: \left(-\frac{\d_o}{2} , \frac{\d_o}{2}\right)^3  \to \bR^3$, such that (a) each point of $\Im(F^{(\d_o, y)})$   is joined to   $y$ through a $\cD^{(\l)}$-path  and (b)  $\Im(F^{(\d_o, y)})$  contains the ball $B_{\mathbf r_o}(y)$. Since  ${\mathbf r_o}$  
  is independent on $y$,  a standard open and closed argument yields that  any  pair $x, x' \in \overline{B_2(0)}$ is  joined by a $\cD^{(\l)}$-path. Combining this  with the fact that  the $X^{(\l)}_i$  are equal to  the $\cC^\infty$ vector fields $X^{(\l = 0)}_i$    on $\bR^3 \setminus B_2(0)$ (that is, on a region  on which $\cD^{(\l)}$ is  a $\cC^\infty$ regular distribution  to which Theorem \ref{Chow} applies), we conclude  that {\it any} pair of   points   $x, x' \in \bR^3$  is  joined by a $\cD^{(\l)}$-path  if   $\l \in (- \d, \d)$,  as claimed. 
 \end{example}
This example illustrates how Corollary \ref{cor43} can be exploited to derive sufficient conditions  for the  existence of  piecewise regular curves,  joining prescribed pairs of points and tangent to  families $\cD = \{\cD_x\}$ of  subspaces $\cD_x \subset T_x \cD$, $x \in \cM$, which  do not  satisfy  the  classical  hypothesis of being generated by  $\cC^\infty$ vector fields.  Results of this kind   are  useful    tools for investigations in Control Theory.  We also remark  that  the structure of the proof of  Theorem \ref{Chow}  yielded the  results in \cite{GSZ1, GSZ2}, where a very similar line of argument  has been  used to get  new local controllability criteria for non-linear real analytic control systems. This suggests  further uses  of Corollary \ref{cor43}  (or, more precisely,  of  appropriate modifications), as for instance  to obtain improvements of  the new    criteria  holding under     regularities assumptions   weaker than real analyticity.

\bigskip
\bigskip
\font\smallsmc = cmcsc8
\font\smalltt = cmtt8
\font\smallit = cmti8
\hbox{\parindent=0pt\parskip=0pt
\vbox{\baselineskip 9.5 pt \hsize=3.1truein
\obeylines
{\smallsmc 
Marta Zoppello
Dipartimento di Scienze Matematiche 
``G. L. Lagrange'' (DISMA)
Politecnico di Torino
Corso Duca degli Abruzzi, 24, 
10129 Torino 
ITALY
}\medskip
{\smallit E-mail}\/: {\smalltt  marta.zoppello@polito.it 
\ 
}
}
\hskip 0.0truecm
\vbox{\baselineskip 9.5 pt \hsize=3.7truein
\obeylines
{\smallsmc
Cristina Giannotti \& Andrea Spiro
Scuola di Scienze e Tecnologie
Universit\`a di Camerino
Via Madonna delle Carceri
I-62032 Camerino (Macerata)
ITALY
\ 
}\medskip
{\smallit E-mail}\/: {\smalltt cristina.giannotti@unicam.it
\smallit E-mail}\/: {\smalltt andrea.spiro@unicam.it}
}
\ 
\ }

\begin{thebibliography}{25}
  
%
%
%
%
%
              
 \bibitem{BGSZ}
F. Bagagiolo, C. Giannotti, A. Spiro and M. Zoppello, {\it Proving the Kalman criterion with flows and orbits}, in preparation.
              
                           
%
   
%
%
%
            
 
            
              
 

 
\bibitem{Ch} W.-L. Chow, {\it \"Uber Systeme von linearen partiellen Differentialgleichungen erster Ordnung}, Math. Ann. {\bf 117} (1939),  98--105.


 
      
      
 
\bibitem{FR}
E. Feleqi and F. Rampazzo,  {\it Iterated {L}ie brackets for nonsmooth vector fields},
NoDEA Nonlinear Differential Equations Appl.   {\bf 24},
(2017), Paper No. 61, pp. 43. 


\bibitem{GSZ1} C. Giannotti, A. Spiro and M. Zoppello, {\it Distributions and controllability problems (I)}, preprint posted in arxiv 2401.07555 (2024). 

\bibitem{GSZ2} C. Giannotti, A. Spiro and M. Zoppello, {\it Distributions and controllability problems (II)}, preprint posted in arxiv  2401.07560 (2024). 


\bibitem{HH}
J. H. Hubbard and B. B. Hubbard, Vector calculus, linear algebra, and differential forms. 
    {\it Prentice Hall, Inc., Upper Saddle River, NJ}, 1999.  
     
     
              
%
\bibitem{Na} T. Nagano,
{\it Linear differential systems with singularities and an
              application to transitive {L}ie algebras},
J. Math. Soc. Japan  {\bf 18},
(1966),
  398--404.
 \bibitem{Ph} G.\ M.\  Phillips, Interpolation and approximation by polynomials, 
 {\it Springer-Verlag, New York}, 2003.
 
 \bibitem{Ra}
P.K. Rashevski\u\i, {\it About connecting two points of complete non-holonomic space by admissible curve} (in Russian), Uch. Zapiski Ped. Inst. K. {\bf 
2} 
(1938), 83--94.
 
 
%
  
  \bibitem{Su} H.\ J.\ Sussmann, {\it Orbits of families of vector fields and integrability of distributions}, Trans. Amer. Math. Soc. {\bf 180}
(1973), 171--188.


\bibitem{SJ} H.\ J.\ Sussmann and V.\ Jurdjevic,
    {\it Controllability of nonlinear systems},  J. Differential Equations  {\bf 12},
(1972), 95--116.


\bibitem{Wa} F.\ W.\ Warner, Foundations of Differentiable Manifolds and Lie groups, {\it Springer-Verlag, New York}, 1983.

 
%
%
%
%

\end{thebibliography}
\end{document}